\documentclass[11pt]{article}
\usepackage{latexsym,amsmath,color,amsthm,amssymb,epsfig,graphicx,mathrsfs}
\usepackage{graphicx}
\input labelfig.tex

\def\qed{\hfill \ifhmode\unskip\nobreak\fi\quad\ifmmode\Box\else$\Box$\fi\\ }
\def\ex{{\rm ex}}

\newtheorem{thrm}{Theorem}[section]
\newtheorem{lemm}[thrm]{Lemma}
\newtheorem{propo}[thrm]{Proposition}
\newtheorem{coro}[thrm]{Corollary}
\newtheorem{defi}[thrm]{Definition}
\newtheorem{conjecture}[thrm]{Conjecture}

\newcommand{\thm}{\begin{thrm}}
\newcommand{\xthm}{\end{thrm}}
\newcommand{\lem}{\begin{lemm}}
\newcommand{\xlem}{\end{lemm}}
\newcommand{\prf}{\begin{proof}}
\newcommand{\xprf}{\end{proof}}
\newcommand{\prop}{\begin{propo}}
\newcommand{\xprop}{\end{propo}}
\newcommand{\cor}{\begin{coro}}
\newcommand{\xcor}{\end{coro}}
\newcommand{\defn}{\begin{defi}}
\newcommand{\xdefn}{\end{defi}}
\newcommand{\conj}{\begin{conjecture}}
\newcommand{\xconj}{\end{conjecture}}

\renewcommand{\phi}{\varphi}

\newtheorem{claim}{Claim}

\newcommand{\bv}{{\bf v}}
\newcommand{\bV}{{\bf V}}

\begin{document}

\title{\vspace{-0.5in}Tur\'an problems and shadows II: trees}

\author{
{\large{Alexandr Kostochka}}\thanks{
\footnotesize {University of Illinois at Urbana--Champaign, Urbana, IL 61801
 and Sobolev Institute of Mathematics, Novosibirsk 630090, Russia. E-mail: \texttt {kostochk@math.uiuc.edu}.
 Research of this author
is supported in part by NSF grant  DMS-1266016
and  by Grant NShÐ1939.2014.1 of the President of
Russia for Leading Scientific Schools.
}}
\and
{\large{Dhruv Mubayi}}\thanks{
\footnotesize {Department of Mathematics, Statistics, and Computer Science, University of Illinois at Chicago, Chicago, IL 60607.
E-mail:  \texttt{mubayi@uic.edu.}
Research partially supported by NSF grants DMS-0969092 and DMS-1300138.}}
\and{\large{Jacques Verstra\"ete}}\thanks{Department of Mathematics, University of California at San Diego, 9500
Gilman Drive, La Jolla, California 92093-0112, USA. E-mail: {\tt jverstra@math.ucsd.edu.} Research supported by NSF Grant DMS-1101489. }}
\maketitle

\vspace{-0.3in}

\begin{abstract}

The expansion $G^+$ of a graph $G$ is the 3-uniform hypergraph obtained from $G$ by enlarging each edge of $G$ with a vertex disjoint from $V(G)$ such that distinct edges are enlarged by distinct vertices. Let $\ex_r(n,F)$ denote the maximum number
of edges in an $r$-uniform hypergraph with $n$ vertices not containing any copy of $F$.  The authors~\cite{KMV} recently determined $\ex_3(n,G^+)$ more generally, namely when $G$ is a path or cycle, thus settling  conjectures of F\"uredi-Jiang~\cite{FJ} (for cycles) and F\"uredi-Jiang-Seiver~\cite{FJS} (for paths).

Here we continue this project by determining the asymptotics for $\ex_3(n,G^+)$   when $G$ is any fixed forest. This settles a conjecture of F\"uredi~\cite{Furedi}.
Using our methods, we also show that
for any graph $G$, either $\ex_3(n,G^{+}) \leq \left(\frac{1}{2} + o(1)\right)n^2$ or $\ex_3(n,G^{+}) \geq (1 + o(1))n^2,$ thereby exhibiting a jump for the Tur\'an number of expansions. 

\end{abstract}

\section{Introduction}

An $r$-uniform hypergraph $F$, or simply {\em $r$-graph}, is a family of $r$-element subsets of a finite set. We associate an $r$-graph $F$ with its edge set and call its vertex set $V(F)$. Given a set of $r$-graphs $\mathcal{F}$, let $\ex_r(n,\mathcal{F})$ denote the maximum number of edges in an $r$-graph on $n$ vertices that does not contain any $r$-graph from $\mathcal{F}$.  When $\mathcal{F}=\{F\}$ we write $\ex_r(n,F)$.   Often we will omit the subscript $r$ in this notation if it is obvious from context, and this paper deals exclusively with the case $r=3$.
Let $G$ be a graph, and for each edge $e \in G$ let $X_e$ be a set of $r - 2$ vertices so that $X_e \cap V(G) = \emptyset$ and $X_e \cap X_f = \emptyset$ when $e \neq f$.
The $r$-uniform {\em expansion} $G^+$ of a graph $G$ is the $r$-graph $G^+ = \{e \cup X_e : e \in G\}$.

\medskip

Expansions of graphs include many important hypergraphs whose extremal functions have been investigated, for instance when $G$ is a triangle
and more generally a clique~\cite{Furedi,FJ,FJS,MV,KMV,M,P2}. Even the simplest case of the expansion of a path with two
edges is non-trivial, in this case we are asking for two hyperedges intersecting in exactly one point. Here the extremal function was determined by
Frankl~\cite{Frankl}, answering a conjecture of Erd\H{o}s and S\'{o}s. 
 If a graph is not $r$-colorable then its $r$-uniform expansion $G^+$ is not $r$-partite, so $\ex_r(n, G^+)=\Omega(n^r)$. We focus on $\ex_r(n,G^+)$ when $G$ is $r$-partite, where a well-known result of Erd\H os yields $\ex(n, G^+) = O(n^{r-\epsilon_G})$ for some $\epsilon_G>0$.

The authors~\cite{KMV} had previously determined $\ex_3(n,G^+)$ exactly (for large $n$) when $G$ is a path or cycle of fixed length $k\ge 3$, thereby answering questions of F\"uredi-Jiang-Seiver~\cite{FJS} and F\"uredi-Jiang~\cite{FJ}.

\subsection{Results}

  A set of vertices in a hypergraph $F$ containing exactly one vertex from every edge of  $F$ is called a {\em crosscut},
following Frankl and F\"{u}redi~\cite{FF}. Let $\sigma(F)$ be the minimum size of a crosscut of $F$ if it exists, i.e.,
\[ \sigma(F) := \min\{|X| : X \subset V(F), \forall e \in F, |e \cap X| = 1\}\]
 if such an $X$ exists.

Since the $r$-graph on $n$ vertices consisting of all edges containing exactly one vertex from a fixed subset of size $\sigma(F) - 1$
does not contain $F$, we have
\begin{equation} \label{triviallower}\ex_r(n,F) \geq (\sigma(F) - 1){n - \sigma(F) + 1 \choose r - 1} \sim (\sigma(F) - 1+o(1)){n \choose r - 1}.\end{equation}
An intriguing open question is when asymptotic equality holds above and this is one of our  motivations for this project. Indeed, it appears that the parameter $\sigma(F)$ often plays a crucial role in determining the extremal function for $F$. The value of $\ex_3(n,G^+)$ was determined precisely by the authors~\cite{KMV} when  $G$ is a  path or cycle.  F\"{u}redi~\cite{Furedi} determined the asymptotics  when $G$ is a forest and $r \geq 4$, by showing that
$\ex_r(n,G^+) = (\sigma(G^+) - 1 +o(1)){n \choose r-1}$. F\"uredi's proof involved extensive use of the delta system method but the method does not work for $r=3$.
Determining $\ex_r(n,G^+)$ seems to get harder as $r$ gets smaller, for example, when $r=2$ and $G$ is a tree, it becomes the Erd\H os-S\'os Conjecture~\cite{ES}.
  F\"uredi conjectured~\cite{Furedi} that $\ex_3(n,G^+) \sim (\sigma(G^+) - 1){n \choose 2}$ when $G$ is a forest, and our main result verifies this conjecture:

\thm \label{main} {\bf (Main Result)}
Let $G$ be a  forest. Then
\[ \ex_3(n,G^+) \sim (\sigma(G^+) - 1){n \choose 2}.\]
\xthm







Our next result concerns $\ex_3(n, G^+)$ for any graph $G$ with $\sigma(G^+) = 2$.

\thm \label{main2}
For every fixed graph $G$ with $\sigma(G^+) = 2$,
 $$\ex_3(n,G^+) \sim {n \choose 2}.$$
\xthm

A straightforward consequence of Theorem \ref{main2} is that for any graph $G$, we have either
$$\ex_3(n,G^+) \leq \left(\frac{1}{2} + o(1)\right)n^2 \quad \mbox{or} \quad \ex_3(n,G^+) \geq (1 + o(1))n^2.$$

\medskip

This paper is organized as follows: in Section \ref{background} we  prove some preliminary lemmas.
In Section \ref{canonicalRamsey}, we give a bipartite version of the canonical Ramsey theorem of Erd\H{o}s and Rado~\cite{ER},
which is one of the main tools for Theorem \ref{main}. We prove Theorem~\ref{main} in Section~\ref{mainsection} and
Theorem \ref{main2} in Section \ref{ktt}. 

\bigskip

{\bf Notation and terminology.} A $3$-graph is called a {\em triple system}. The edges will be written
as unordered lists, for instance, $xyz$ represents $\{x,y,z\}$.
For a set $X$ of vertices of a hypergraph $H$, let $H - X = \{e \in H : e \cap X = \emptyset\}$.  If $X = \{x\}$, then we write
$H - x$ instead of $H - X$.
The {\em codegree} of a pair $\{x,y\}$ of vertices in $H$ is $d_H(x,y) = |\{e \in H : S \subset e\}|$
and for a set $S$ of vertices, $N_H(S)=\{x \in V(H) : S \cup \{x\} \in H\}$ so that $|N_H(S)|=d_H(S)$ when $|S| = 2$.
The {\em shadow} of $H$ is the graph $\partial H = \{xy : \exists e \in H, \{x,y\} \subset e\}$.
The edges of $\partial H$ will be called the {\em sub-edges} of $H$.
A triple system is {\em linear} if every pair of its edges intersect in at most one point.
For an edge $e$ in a triple system $H$, let $\delta_H(e)$ and $\triangle_H(e)$ respectively denote the smallest and largest codegree among the three pairs in $e$.

\section{Expansions in sparse hypergraphs}\label{background}

In this section we state and prove a basic  result about hypergraphs that  generalizes  the fact that a graph with average degree $d$ contains a subgraph of minimum degree at least $d/2$.


\defn
A triple system  $H$ is {\em $d$-full} if every sub-edge of $H$ has codegree at least $d$.
\xdefn

Thus $H$ is $d$-full is equivalent to saying that the minimum non-zero codegree in $H$ is at least $d$.
The following lemma extends the well known fact that any graph $G$ has a subgraph of minimum degree at least $d+1$ with at least $|G| - d|V(G)|$ edges.

\lem \label{fullsub}
For  $d \geq 1$, every $n$-vertex triple system $H$ has a $(d + 1)$-full subgraph $F$ with
\[ |F| \geq |H| - d|\partial H|.\]
\xlem

\prf
A {\em $d$-sparse sequence} is a maximal sequence $e_1,e_2,\dots,e_m \in \partial H$ such that $d_H(e_1) \leq d$, and for all $i > 1$, $e_i$
 is contained in at most $d$ edges of $H$ which contain none of $e_1,e_2,\dots,e_{i - 1}$.
The $3$-graph $F$ obtained by deleting all edges of $H$ containing at least one of the $e_i$ is $(d + 1)$-full. Since
a $d$-sparse sequence has length at most $|\partial H|$, we have $|F| \geq |H| - d|\partial H|$.
\xprf

\section{Colors, lists, and a canonical Ramsey theorem}\label{canonicalRamsey}

One of our main new tools is to use the canonical Ramsey theorem of Erd\H{o}s and Rado~\cite{ER}.
We need a bipartite version of this classical result.

\defn Let $F$ be a bipartite graph with parts $X$ and $Y$ and an edge-coloring
$\chi$. Then 

\begin{tabular}{lp{5.5in}}
{\rm 1.} &  $\chi$ is {\em $X$-canonical} if for each $x \in X$, all edges of $F$ on  $x $ have the same color and edges on different vertices in $X$ have different colors \\
{\rm 2.} &  $\chi$ is {\em canonical} if $\chi$ is $X$-canonical or $Y$-canonical \\
{\rm 3.} &  $\chi$ is {\em rainbow} if the colors of all the edges of $F$ are different and \\
{\rm 4.} &  $\chi$ is {\em monochromatic} if the colors of all the edges of $F$ are the same.
\end{tabular}
\xdefn

Recall that a {\em sunflower} or {\em $\Delta$-system} is  a collection of sets such that the intersection of any two of them is equal to the intersection of all of them.
A key result on sunflowers is the Erd\H{o}s-Rado Sunflower Lemma~\cite{sunflower}:

\lem {\bf (Erd\H{o}s-Rado Sunflower Lemma)} \label{sunlemma}
If $F$ is a collection of sets of size at most $k$ and $|F|\geq k!(s-1)^{k}$,
then $F$ contains a sunflower with $s$ sets.
\xlem

If $\chi$ is an edge-coloring of a graph $F$ and $G \subset F$, then $\chi\vert_G$ denotes the edge-coloring of $G$ obtained by restricting
$\chi$ to the edge-set of $G$. A bipartite version of the canonical Ramsey theorem is as follows:

\thm \label{canonical} For each $s>0$ there exists $t>0$ such that for any  edge-coloring $\chi $ of $G = K_{t,t}$,
 there exists  $K_{s,s} \subset G$ such that $\chi \vert_{K_{s,s}} $ is monochromatic or rainbow or canonical.
\xthm

\prf Let $X$ and $Y$ be the parts of $G$ and let $S = \{y_1,y_2,\dots,y_{2s^2}\} \subset Y$. Let $W$ be the set of
vertices $x \in X$ contained in $s$ edges of the same color connecting $x$ with $S$. If $|W| > m := s^2{2s^2 \choose s}$, then
there is a set $Y' \subset S$ of size $s$ and a set $X' \subset W$ of size $s^2$ such that for every
$x \in X'$, the edges $xy$ with $y \in Y'$ all have the same color. In this case we
recover either a monochromatic $K_{s,s}$ or an $X'$-canonical $K_{s,s}$. Now suppose $|W| \leq m$.
For $x \in X_0 := X \backslash W$, let $C(x)$ be a set of $2s$ distinct colors on edges between $x$ and $S$.
By the Sunflower Lemma, Lemma \ref{sunlemma}, if $|X_0| > (2s)! (s!m)^{2s}$, then there exists $X_1 \subset X_0$
such that $\{C(x) : x \in X_1\}$ is a $\Delta$-system of size $s!m$. Let $C$ be the core of this $\Delta$-system.
First suppose $|C| \geq s$. Then we have a set $X_2$ of at least $s \cdot s!$ vertices in $X_1$ which each sends
 $s$ edges with colors from $C$ into a fixed subset $Y_3$ of $S$ of size $s$.
 This implies that for some set $X_3 \subset X_2$ of size $s$, the $K_{s,s}$ between $X_3$ and $Y_3$ is
$Y_3$-canonical. Finally, suppose $|C| < s$. Pick $C'(x) \subset C(x) \backslash C$ of size $s$ for  $x \in X_0$.
Since $|X_0| = s!m > s{2s^2 \choose s}$, we find a set $Y^* \subset S$ of size $s$ as well as a set $X^* \subset X$ of $s$ vertices
$x \in X_0$ such that the edges between $x$ and $Y^*$ have colors from $C'(x)$. Since the sets $C'(x)$ are disjoint,
 this is a rainbow copy of $K_{s,s}$. \xprf


We now link this to the context of hypergraphs via the following definition.

\defn Let $H$ be a $3$-graph. For $G \subset \partial H$ and $e \in G$, let
\[ L_G(e) = N_H(e)\setminus V(G).\]
The set $L_G(e)$ is called the {\em list} of $e$ and the elements of $L_G(e)$ are called {\em colors}.
\xdefn

Let $L_G = \bigcup_{e \in G} L_G(e)$ -- this is the set of colors in the lists of edges of $G$.

\defn A {\em list edge coloring of $G$} is a map $\chi : G \rightarrow L_G$ with $\chi(e) \in L_G(e)$
for all $e \in G$. List-edge-colorings $\chi_1,\chi_2 : G \rightarrow L_G$ are {\em disjoint} if
$\chi_1(e) \neq \chi_2(f)$ for all $e,f \in G$.
\xdefn

If $\chi$ is an injection -- the coloring is rainbow -- then clearly $G^{+} \subset H$. We require one more definition:

\defn Let $H$ be a $3$-graph and $m \in \mathbb N$. An {\em $m$-multicoloring} of $G \subset \partial H$ is a family
of list-edge-colorings $\chi_1,\chi_2,\dots,\chi_m : G \rightarrow L_G$ such that $\chi_i(e) \neq \chi_j(e)$ for every $e \in G$
and $i\neq j$.
\xdefn

A necessary and sufficient condition for the existence of
an $m$-multicoloring of $G$ is that all edges of $G$ have codegree at least $m$ in $H$.
We stress here that the definitions are all with respect to the fixed host $3$-graph $H$. The following result will be key to the proofs of Theorem \ref{main} and Theorem \ref{main2}.

\thm \label{canon} Let $m,s \in \mathbb N$, let $H$ be a $3$-graph, and let $G = K_{t,t} \subset \partial H$.
Suppose $G$ has an $m$-multicoloring. If $t$ is large enough, then there exists $F = K_{s,s} \subset G$ such that $F$ has
either a rainbow list-edge-coloring or an $m$-multicoloring such that the colorings are pairwise disjoint, and each
coloring is monochromatic or canonical.
\xthm

\prf Set $s=t_m/m^2$ and $t_m < t_{m-1} < \cdots < t_1 < t_0=t$ where Theorem \ref{canonical} with input $t_i$ has output $t_{i-1}$.
Pick a  color $c_1(e)$ on each edge $e\in G$ and apply Theorem \ref{canonical} to $G$.   We obtain a
rainbow, monochromatic or canonical subgraph $G_1$ of $G$ where $G_1=K_{t_1, t_1}$. If it is rainbow, then we are done, so assume it is monochromatic or canonical. For every $e\in G_1$, remove $c_1(e)$ from its list. Now pick another color on each edge of $G_1$ and repeat. We obtain subgraphs
$G_m \subset G_{m-1} \subset \cdots \subset G_1$ such that each $G_i$ is monochromatic or canonical where $G_i=K_{t_i, t_i}$ has parts $X_i, Y_i$.
In particular, each coloring of the $m$-multicoloring of $G$ restricted to $G_m$ is monochromatic or canonical.
\smallskip

Let us assume that we have $a$ monochromatic colorings, $b$ $X_m$-canonical colorings, and $c$ $Y_m$-canonical colorings of $G_m$ where $a+b+c=m$.
It suffices to ensure that these colorings are pairwise disjoint.
A  color $\chi(xy)$ in an $X_i$-canonical coloring of  $G_i$ cannot appear in a $Y_{i'}$-canonical coloring of $G_{i'}$ for $i'>i$ as $\chi(xy)$ was deleted from all edges incident to $x$ when forming $G_{i+1}$.  A similar statement holds with $X$ and $Y$ interchanged, so every $X_m$-canonical coloring of $G_m$ is disjoint from every  $Y_m$-canonical coloring of $G_m$.  The same argument shows that no color in a monochromatic coloring appears in a canonical coloring. It suffices to show that colors on different $X_m$-canonical colorings are disjoint (and the same for $Y_m$-canonical).
\smallskip

Let the $b$  $X_m$-canonical colorings be
$\chi_1,...,\chi_b$.
Construct an auxiliary graph $K$ with $V(K)=X_m$ where $xx'\in K$ if there exist $i\ne i'$ and a  color
$\alpha$ that is canonical for  $x$ in $\chi_i$ and  canonical for $x'$ in $\chi_{i'}$.
Then $\Delta(K)\leq b(b-1)$, so $K$ has an independent set of size
$s=t^m/m^2\le |X_m|/b^2$. Let us restrict $X_m$ to this independent set. We repeat this procedure for $Y_m$ and finally obtain a subgraph $F=K_{s,s}$ with an $m$-multicoloring that satisfies the requirement of the theorem.
\xprf

\section{Cleaning lemmas}\label{local}

The  lemmas in this section will allow us to find for an appropriate
triple system $H$ a large dense graph $G \subset \partial H$ that possesses an $m$-multicoloring with the colors outside of $V(G)$.
Using such substructures, we will  embed expansions of graphs into $H$.

\lem \label{kttlemma}
Let $m,t \in \mathbb N$, $\delta \in \mathbb R_+$ and $H$ be an $n$-vertex triple system.
Suppose that $F \subset \partial H$ and for each $f \in F$ let $S_f \subset V(H)\setminus f$  with $|S_f|=m$.
 If $|F| \geq \delta n^{2}$ and $n$ is large enough, then there exists $K \subset F$ such that $K \cong K_{t,t}$
 and $S_f \cap V(K)=\emptyset$ for each $f \in K$.
\xlem

\prf Let $T$ be a random subset of $V(H)$ obtained by picking each vertex independently with probability $p=1/2$. Let
$G = \{f \in F : f \subset T, S_f \cap T=\emptyset\}.$
Then
$$\mathbb E(|G|) \ge |F| p^{2}(1-p)^{m}\ge \frac{\delta}{2^{m+2}}n^{2}.$$
So there is a $T \subset V(H)$ with $|G|$ at least this large.
If $n$ is large enough, then the K\"{o}vari-S\'{o}s-Tur\'{a}n Theorem implies that there exists a complete bipartite graph $K \subset G \subset F$
with parts of size $t$. Due to the definition of $G$, the subgraph $K$ satisfies the requirements of the lemma.
\xprf

\lem \label{trimlemma}
Let $A_1,\ldots,A_m$ be disjoint subsets of a set $V$ and $a_1,\ldots,a_m$ be distinct elements of $V$.
Then there are $\lceil\frac{m}{3}\rceil$ pairwise disjoint sets of the kind $A_i+a_i$.
\xlem

\prf Note that the statement of the lemma allows $a_i\in A_i$. Since all $a_1,\ldots,a_m$ are distinct, if $(A_i+a_i)\cap (A_j+a_j)\neq \emptyset$, then $a_i\in A_j$, or $a_j\in A_i$, or both.
Let $F$ be the digraph with vertex set $\{A_1+a_1,\ldots,A_m+a_m\}$ and $A_iA_j\in F$ if $a_i\in A_j$ and $i\neq j$. Since the outdegree of every
vertex in $F$ is at most $1$, $F$ is $3$-colorable and thus has an independent $I$ of size $\lceil\frac{m}{3}\rceil$.
By definition, the members of $I$ are  pairwise disjoint.
\xprf

\section{Trees and crosscuts}
In this section we produce a structural decomposition of a tree $T$ that will be used later to embed $T^+$ in a hypergraph.  We will also prove some basic lemmas about this decomposition.

Let $G$ be a graph and consider a minimum crosscut $X$ of $G^+$. For an edge $e \in G$, let $v_e$ be the unique vertex such that $e \cup v_e \in G^+$; say that $v_e$ is the {\em enlargement} of $e$.  Partition $X$ into $I \cup J$ where $J$ comprises the vertices of $X$ that are used for enlargement of the edges of $G$.  Let $R \subset G$ be the set of edges $e$ such that $v_e \in J$. Then  $I \subset V(G)$ is an independent set in $G$ and $R \subset G-I$.  Furthermore,
$\sigma(G^+)=|X|=|I|+|R|$. On the other hand for any independent set $I \subset V(G)$ and subgraph  $R \subset G-I$, such that every edge of $G-I$ is in $R$, we obtain a crosscut $X= I \bigcup \cup _{e \in R} \{v_e\}$ of $G^+$. Consequently,
$$\sigma(G^+)=\min \{|I|+|G-I|: I \subset V(G) \hbox{ is an independent set}\}.$$
In the ensuing proof, it is more convenient to work with the pair $(I, R)$ rather than a crosscut $X$ of $G^+$.
\medskip

\defn
A {\em crosscut pair} of a graph $G$ is a pair $(I,R)$ where \\
\begin{tabular}{lp{5in}}
$\bullet$ & $I \subset V(G)$ is an independent set, \\
$\bullet$ & $R = \{e \in G : e \cap I= \emptyset\}$.
\end{tabular}
\newline The crosscut pair $(I, R)$ is {\em optimal} if $|I|+|R|=\sigma(G^+)$.
\xdefn
 Given a crosscut pair $(I, R)$, let
$$ L = \{v \in V(G) \backslash (V(R) \cup I) : d_{G}(v)= 1\} \quad and \quad
D =\{v \in V(G) \backslash (V(R) \cup I) : d_{G}(v)> 1\}$$
so that $D=V(G) \setminus (V(R) \cup I \cup L)$.

\medskip

\begin{center}
\centerline{\includegraphics[width=5.8in]{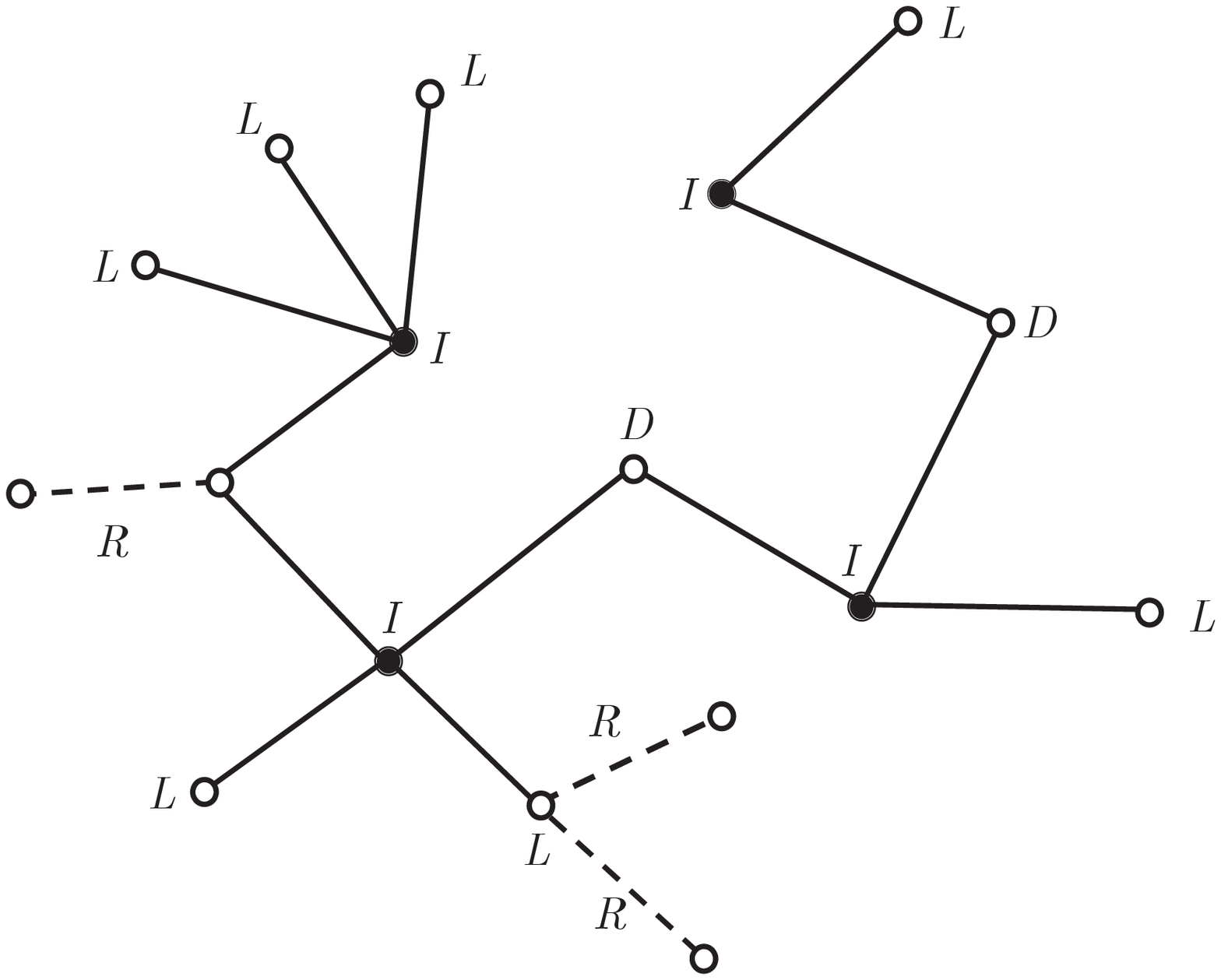}}

\vspace{-3.4in}

Figure 1 : Optimal crosscut decomposition $\sigma(G^+) = 7 = |I| + |R|$
\end{center}

\medskip

\lem \label{jl1}
Let $T$ be a  tree  with $\sigma(T^+)=\ell+1>0$.
Consider an optimal crosscut pair $(I,R)$ of $T$ that maximizes $|I|$.
Then  $(a)$ $|R|\leq \ell/2$ and $(b)$ no pendant edge of $T$ belongs to $R$.
\xlem
{\bf Proof.} Suppose $yz \in R$ is a pendant edge of $T$ with $d_T(z) = 1$. Then $z$ is not adjacent to any vertex of $I$.
The crosscut pair $(I',R')$ where $I' = I \cup \{z\}$ and $R' = R \backslash \{yz\}$
contradicts the maximality of $I$, proving (b). Since $T$ has no cycles, counting the edges
induced by $I \cup V(R)$, we have
\[ |V(R)| + |R| \leq |I| + |V(R)| - 1\]
This yields $|I|\geq 1+|R|$. Since $(I,R)$ is optimal, $|I| + |R| = \ell + 1$ and this gives (a).\qed

\medskip

\lem \label{forest-tree}
Let $F$ be a $k$-vertex forest. Then there exists a $k$-vertex tree $T \supset F$ with $\sigma(T^+)= \sigma(F^+)$.
\xlem

{\bf Proof.} Let $F$ have components $T_1,T_2, \ldots, T_s$. For each $j=1, \ldots, s$, let $(I_j, R_j)$ be an optimal crosscut pair of $T_j$ and with
 $|I_j|$ a maximum. If $I_j = \emptyset$, then any pendant edge of $T_j$ is in $R$, contradicting Lemma \ref{jl1}.(b). Therefore
 $I_j \neq \emptyset$ for $j = 1,2,\dots,s$. If $I = \bigcup_{j=1}^s I_j$ and $R = \bigcup_{j=1}^s R_j$, then clearly $(I,R)$
is an optimal crosscut pair of $G$, and $\sigma(G^+) = \sum_{j = 1}^s \sigma(T_j^+)$. For each $j : 1\le j \le s-1$, let us add an edge between $T_j$ and $T_{j+1}$ as follows: pick $u \in V(T_j)-I_j$ and $v \in I_{j+1}$ and add the edge $uv$.  This results in a tree $T \supset F$ with $\sigma(T^+) =\sigma(F^+)$. \qed

\medskip

\defn
Let $T$ be a tree with parts $P$ and $Q$, with $|P|\leq |Q|$. Then
\[ \lambda(T) = \left\{\begin{array}{ll}
|P| - 1 & \mbox{ if some leaf of }T\mbox{ is in }P\\
|P| & \mbox{ otherwise}
\end{array}\right.\]
If $F$ is a forest with components $S_1,\ldots,S_h$, then we define $\lambda(F)=\sum_{i=1}^h\lambda(S_i)$.
\xdefn

\lem \label{lambda-lemma}
For every forest $F$, $\lambda(F) \leq |F|/2$.
\xlem

{\bf Proof.} It is known that if $T$ is a tree with parts $P$ and $Q$ and $|P|=|Q|$, then each of $P$ and $Q$ contains a leaf. This shows $\lambda(T) \leq |T|/2$ for
every tree $T$, and applying this to the components of $F$, we obtain the lemma. \qed

\lem \label{treelem}
Let $T$ be a tree with an optimal crosscut pair $(I,R)$ and $\sigma(T^+) = \ell + 1 > 0$. 
Then
 \begin{equation}\label{root}
 d_{T}(r)\leq  \ell-\lambda 
  \quad\mbox{for each 
 $r\in V(R)$.}
 \end{equation}
 \xlem

 {\bf Proof.} Suppose that $R$ consists of $h$ components $R_1, \ldots, R_h$ and $r\in V(R_1)$. The second end
 of every edge $rv\in G-R$ must be in $I$. Also every vertex of $R$ has a neighbor in $I$,
 because otherwise we could move the vertex into $I$ and obtain a crosscut pair of the same size and larger $|I|$.
 Moreover,  as $T$ is acyclic, every $V(R_j)$ has at least $|V(R_j)|=1+|R_j|$ neighbors in $I$ and the neighborhoods of sets $R_j$
 in $I$ form a hypergraph linear forest. Thus
 $r$ is not adjacent
 to at least $\sum_{j=1}^h|R_j|=\ell_2$ vertices in $I$.
 By definition, $r$ is not  adjacent to at least $\lambda(R_1)$
 vertices in $V(R_1)$ (the smaller partite set of $V(R_1)$). So, since $|R|-\lambda\geq |R_1|-\lambda(R_1)$,
 \begin{equation}\label{s22}
 d_{T}(r)\leq (|I|-\ell_2)+(|V(R_1)|-\lambda(R_1))\leq \ell_1-\ell_2+(1+|R|-\lambda)=\ell_1+1-\lambda.
 \end{equation}
The last expression is at most  $\ell-\lambda$
unless $\ell_2=1$.  Suppose  $\ell_2=1$. Then $h=1$, $\lambda=0$, $r$ has exactly one neighbor in $V(R)$, and instead of~\eqref{s22},
we have
$$d_{T}(r)\leq (|I|-\ell_2)+1=\ell_1=\ell-\lambda.$$
 So,~\eqref{root} holds in this case, as well.
\qed

\section{Proof of Theorem \ref{main}} \label{mainsection}

Let $G$ be a forest with $k$ vertices, and $\ell = \sigma(G^+) - 1$. We are to show $\ex_3(n,G^+) \leq (\ell + o(1)){n \choose 2}$. Let $H$
be a triple system on $n$ vertices with $|H| > (\ell + \epsilon){n \choose 2}$ where $\epsilon > 0$.
By Lemma \ref{forest-tree}, $G^+ \subset T^+$ for some tree $T$ with $k$ vertices and $\sigma(T^+) = \sigma(G^+)$, so it is enough to show $T^+ \subset H$.
for $n > n_0(\epsilon,k)$. Suppose for a contradiction that $T^+ \not \subset H$.

\subsection{Finding a rich triple system}
Recall that  $\delta_H(e) = \min_{uv \subset e} d_H(uv)$.
In this section we show how to find hypergraphs $H_3 \subset H_1 \subset H$ such that  $\delta_{H_1}(f)\ge \ell+1$ for every $e\in H_1$,
$H_3$ has quadratically many edges,
and  $\delta_{H}(f)\le 3k$ for all $f\in \partial H_3$.  We will later use $H_1$ and $H_3$ to embed $T^+$.


Let $H_1$ be obtained from $H$ by consecutive deletion of edges having a pair
of codegree at most $\ell$ in the current $3$-graph, so that $\delta_{H_1}(e) \geq \ell + 1$ for all $e \in H_1$.
Let $F_1$ denote the set of deleted edges.
Since we deleted at most $\ell$ edges at each step and the number of steps is at most ${n \choose 2}$, we have
$|F_1|\le \ell n^2/2$ and
$|H_1| = |H_0| - |F_1| \geq (\ell+\epsilon)n^2/2-\ell n^2/2 \ge \epsilon n^2/2.$ By definition,
\begin{equation}\label{0119}
\delta_{H_1}(e) \geq \ell+1
\qquad\mbox{for every $e \in H_1$.}
\end{equation}
Let
\[ H_2 = \{e \in H_1 : \delta_{H}(e) > 3k\} \quad \mbox{ and } \quad  H_3 = \{e \in H_1 : \delta_{H}(e) \leq 3k\},\]
so that $H_1 = H_2 \cup H_3$. Suppose for a contradiction that $|H_2| > 3k^2 n$. If $|\partial H_2| > kn$, then $T \subset \partial H_2$,
and we greedily extend $T$ to $T^+ \subset H$. Otherwise, $|\partial H_2| \leq kn$, in which case
$H_2$ has a  $3k$-full subgraph of size at least $|H_2| - 3k|\partial H_2| > 0$, which clearly contains $T^+$.
This contradiction shows $|H_2| \leq 3k^2 n$, and therefore $|H_3| = |H_1| - |H_2| \geq \epsilon n^2/4$ for large enough $n$.


By the definition of $H_3$,
in each $e \in H_3$ we can fix some $f_e \in {e \choose 2}$
with $d_{H}(f_e)\le 3k$. Let $F=\{f_e: e \in H_3\}$.  Then $|F|\ge |H_3|/3k> \epsilon n^2/12k$.
For each $f \in F$, let $S_f \subset N_{H_1}(f)$ with $|S_f|=\ell+1$.
Applying  Lemma~\ref{kttlemma} to $F \subset \partial
H_3$ we find a copy $K$ of $K_{t,t}$ for large $t$ such that each
edge $f$ of  $K$ is contained in  $\ell+1$ edges $f \cup \{v\} \in H_1$ with $v \in S_f$.

The $\ell+1$ edges $f \cup \{v\} $ with $v \in S_f$ containing   $f$ for every $f\in K$ give an $(\ell+1)$-multicoloring of $K$,
so by Theorem~\ref{canon}
there is $G_0=K_{s,s} \subset K$ ($s$ large) with an $(\ell+1)$-multicoloring $M_1,\ldots,M_{\ell+1}$ such that

\begin{center}
\begin{tabular}{lp{5in}}
$\bullet$ & some $M_i$ is rainbow, or \\
$\bullet$ & the $M_i$'s are pairwise disjoint and each $M_i$ is  either monochromatic or canonical.
\end{tabular}
\end{center}

 Let $X$ and $Y$ be the partite sets of $G_0$ and
 \[ Z=\bigcup_{x\in X, y \in Y} N_{H_1}(xy) - V(G_0).\]

We will often think of $M_i$ as a 3-graph comprising the edges $xyw$ where $x \in X, y \in Y$ and $w$ is the color of $xy$.

\subsection{Canonical colorings and embeddings}
In this section we prove a series of claims using Theorem~\ref{canon}
that allow us to embed $T^+$ within $H_1$ in certain situations.

\begin{claim}\label{cla1} No $M_i$ is  rainbow.
\end{claim}

{\em Proof:} Suppose $M_1$ is rainbow.
Since $s>3k$, there is an embedding $\psi(T)$ of $T$ into $G_0$. Since $M_1$ is rainbow, its edges containing the edges
of $\psi(T)$ form $T^+ \subset H_1$.\qed

\begin{claim}\label{cla3} If some $M_i$ is  $Y$-canonical then there are no $X$-canonical $M_j$.
\end{claim}
{\em Proof:} Suppose $M_1$ is  $Y$-canonical and $M_2$ is  $X$-canonical.
Then for every $y\in Y$ there is $w(y)$ such that $xyw(y) \in H_1$ for all $x \in X, y \in Y$ and
 for every $x\in X$ there is $u(x)$ such that $xyu(x) \in H_1$ for all $x \in X, y \in Y$.
Let $\hat{T}$
be a directed out-rooted tree obtained from $T$ with any root $v$. We embed it
into $G_0$, and expand each edge as follows: If the image of a directed edge of $\hat{T}$ is
$xy$, then expand it to $xyw(y)$ and if the image is $yx$, then expand it to $yxu(x)$.\qed

Choose an optimal crosscut pair of $T$ with maximum $|I|$. Let $\ell_1=|I|$ and $\ell_2 = |R|$.
By Lemma~\ref{jl1}.(b), the pendant edges of $T$ are not in $R$.

\begin{claim}\label{cla4} At most $\ell_1-1$ of the $M_i$ are monochromatic.
\end{claim}
{\em Proof:} Suppose, without loss of generality, that for $i=1,2,\ldots,\ell_1$, each $M_i$ is monochromatic and $w_i$
is the common vertex of all edges in $M_i$. If $I=\{a_1,\ldots,a_{\ell_1}\}$, then for
$i=1,\ldots,\ell_1$, we place $a_i$ onto $w_i$, and then embed $T-I$ into $G_0$.
Since each of $w_1,\ldots,w_{\ell_1}$ is adjacent in $\partial H_1$ with each vertex of $G_0$,
this yields an embedding of $T$ into $\partial H_1$.
Next we extend the $\ell_2$ edges of $R$ using for each of them an edge from one of the $\ell_2$ sets
$M_{\ell_1+1},\ldots,M_{\ell_1+\ell_2}$ (one edge from each set). Every other edge of the embedded $T$ is
incident with one of $w_i$. If such an edge has the form $w_ix$ (respectively, $w_iy$) then we take any unused
$y\in Y$ (respectively, $x\in X$) and extend it to $\{w_i,x,y\}$.\qed

\begin{claim}\label{cla4a} $R\neq\emptyset$.
\end{claim}
{\em Proof:} Suppose $R=\emptyset$ and $U,U'$ are partite sets of $T$.
Then all vertices of $I$ are in the same partite set, say $U$, of $T$ (in fact, $I=U$ as $I$ covers all edges of $T$).
By Claims~\ref{cla1},~\ref{cla3} and~\ref{cla4} and symmetry,
we may assume that $M_1$ is $Y$-canonical.
For every $y\in Y$ there is $w(y)$ such that $xyw(y) \in H_3$ for all $x\in X, y \in Y$.
Place the vertices of $T$ into $X\cup Y$ so that
$U\subset X$ and $U' \subset Y$. Since $G_0$ is a complete bipartite graph, this yields an embedding of $T$ into $G_0$.
Since $T$ is a tree, $|T|=|U'|+|U|-1=|U'|+\ell$. For every $y\in Y$ which is the image of some $b\in U'$ we expand
one edge $xy$ by adding $w(y)$. For the remaining $\ell$ edges we use edges of $M_2,\ldots,M_{\ell+1}$,
from distinct $M_j$ for distinct edges.
\qed

We recall from the last section the definition of $\lambda(F)$ for a forest $F$.
\begin{claim}\label{cla5} If some $M_i$ is $Y$-canonical, then  at most $\lambda(R)-1$ of the $M_j$ are monochromatic.
\end{claim}
{\em Proof:} Suppose that $M_\ell$ is $Y$-canonical (we may assume this by Claim 3) and suppose, for a contradiction, that for $i=1,\ldots, \lambda(R)$,
each $M_i$ is monochromatic and $w_i$
is the common vertex of all edges in $M_i$.
Also for every $y\in Y$ there is $w(y)$ such that each edge in $M_\ell$ containing $y$ also contains $w(y)$. We embed $R$
into the subgraph of $\partial H_1$ induced by $Y\cup \{w(y)\,:\,y\in Y\}\cup \{w_1,\ldots,w_{\lambda(R)}\}$
as follows. Suppose the components formed by the edges of $R$ are $R_1,\ldots,R_h$ with smaller partite sets $P_1,\ldots,P_h$
and if $P_j$ contains leaves, then $b_j$ is one of them. We choose arbitrary $y_1,\ldots,y_h\in Y$, and for $j$ such that $b_j$
exists, place $b_j$ onto $w(y_j)$ and the neighbor in $R_j$ of $b_j$ onto $y_j$.
Then place the remaining $\lambda(R)$ vertices of $P_1\cup\ldots\cup P_h$ onto vertices
in $\{w_1,\ldots,w_{\lambda(R)}\}$ and the remaining vertices of $V(R)$
(which comprise $\bigcup _i V(R_i)\setminus(P_i \cup N_{R}(b_i))$) onto arbitrary free vertices in $Y$.
Since each $w_jy \in \partial H_1$ for all $y\in Y$ this yields an embedding of $R$ in $\partial H_1$. Next, place the vertices of $D \cup L$ into  new free
vertices of $Y$,  and finally, place all vertices of $I$ onto distinct vertices in $X$.

This gives an embedding of $T$ into $\partial H_1$. We expand it to an embedding of $T^+$ into $H_1$ as follows.
Since   $xyw(y)\in H_1$ for all $x \in X$ and $|X|\ge s$, we can expand the edges of the form $yw(y)$ at the end. Expand all edges of the form $w_iy$ and $w_ix$
 by adding a free vertex from $X$ and $Y$, respectively. This allows us to expand all edges of $T$ except those that contain some vertex of $D$ as an endpoint.
We now focus on these edges which connect $D$ to $I$.

For every $y$ onto which we placed a vertex $a\in D$,
we expand one
edge of the kind $xy$ by adding $w(y)$ and all other such edges
using some $M_j$ (distinct for distinct $ab$). To prove that we have enough free $M_j$ first observe that
the number of edges in $T$ connecting $D$ to $I$ is
$|D|+\ell_1-1$ (because $I$ is an independent set and each edge joins precisely two components). Of these edges, $|D|$ will be expanded by expanding pairs of the form $yw(y)$ as mentioned earlier, so we must only expand $\ell_1-1$ more edges.
 The number of $M_j$ that have already been used is at most $\lambda(R)+1$ and so the number of unused $M_j$
 is at least $\ell+1-(\lambda(R)+1)= \ell-\lambda(R)$.
We finally show $\ell - \lambda(R) \ge \ell_1 - 1$ to complete the embedding.
By Lemma \ref{lambda-lemma}, $\lambda(R) \leq |R|/2$. Since $|I| + |R| = \ell + 1$ and $|R| \leq \ell/2$ from Lemma \ref{jl1},
we have $\ell - |R|/2 \geq \ell_1 - 1$, and therefore
$$\ell-\lambda(R)\ge \ell -\frac{|R|}{2} \ge \ell_1 - 1.$$
This shows that $T^+ \subset H$, a contradiction. \qed

\medskip

Let $\lambda=\lambda(R)$.
By Claims~\ref{cla1},~\ref{cla3},~\ref{cla4} and symmetry,
we may order the colorings so that $M_1,\ldots,M_{p}$ are $Y$-canonical and the remaining are monochromatic.
Furthermore, by Claim~\ref{cla5},
$$p\geq {\ell+2-\lambda}.$$

\subsection{Constructing the digraph $D_g$}
In this section we construct a digraph $D_g$ whose underlying edges lie in $\partial H_1$.  The digraph $D_g$ will be the vertex disjoint union of homomorphic images of directed out-trees each with height at most $\ell_2+1$. The rich structure of $D_g$ encodes edges of $H_1$ and will later be used to embed $T^+$ in $H_1$.

 For each $i\in \{1,\ldots,p\}$
and every $y\in Y$, let $w_i(y)$ be the vertex such that each edge in $M_i$ containing $y$ contains also $w_i(y)$ and let
$W_i=\{w_i(y)\,:\,y\in Y\}$. Also for every $y\in Y$, let $W(y)=\{w_i(y)\,:\,i=1,\ldots,p\}$. Let $Q=\{\alpha_{p+1},\ldots,
\alpha_{\ell+1}\}$ be the set of the  colors used in the monochromatic colorings $M_i$.

 By definition, for each $i\in \{1,\ldots,p\}$, the subgraph of $\partial H_1$ induced by
$X \cup W_i$ contains the complete bipartite graph with partite sets $X$ and $W_i$. By Theorem~\ref{canon},
all $W_i$ are mutually disjoint and disjoint from $Y$. 
By the same theorem, we also have $W_i \cap Q=\emptyset$ where $Q$ is the set of vertices/colors  in the monochromatic colorings
$M_{p+1}, \ldots, M_{\ell+1}$.

{\bf Basic cleaning procedure:} By~\eqref{0119}, for each $x \in X, w_1(y) \in W_1$, we can choose a set $S(x,w_1(y)) \subset N_{H_1}(xw_1(y))$
with $y \in S(x,w_1(y))$ and $|S(x, w_1(y))|=\ell+1$.
 Define the 3-graph
 $$H_1'=\{xwz \in H_1: x \in X, w \in W_1, z \in S(x,w)\}$$ with $V(H_1')=\cup_{e \in H_1'} e$ so that
$$|V(H_1')|\le |X|+|W_1|+ (\ell+1)|X||W_1|< (\ell+3)s^2.$$
Let $F'$ be the complete bipartite graph with parts $X$ and $W_1$ so that $F' \subset \partial H_1'$.
Then $|F'|=|X||W_1|=s^2\ge\delta |V(H_1')|^2$ for  $\delta=1/(\ell+3)$. Since $s$ is large,  we may apply
Lemma~\ref{kttlemma} to $F' \subset \partial H_1'$ to obtain a large complete bipartite subgraph
 $G_{1,1}\subset F'$ such that $S(x,w)\cap V(G_{1,1})=\emptyset$ for all $xw \in G_{1,1}$.
Since $|S(x,w)|\ge \ell+1$ for all $xw$, we can view $G_{1,1}$ as being multicolored with $\ell+1$ colors, with one of the color classes corresponding to the vertices $y$. Moreover, all colors lie outside $V(F')$.
The reason  we need this is to apply Claim 1 below.
This is the {\em basic cleaning procedure}.

By Theorem~\ref{canon}, we obtain subsets
 $X'_{1,1}\subset X, W'_{1,1}\subset W_1$, such that the $(\ell+1)$-multicoloring restricted to $X'_{1,1} \times W'_{1,1}$ comprises rainbow, monochromatic, or canonical colorings. Let $Y'_{1,1}=\{y \in Y: w_1(y) \in W'_{1,1}\}$. None of the colorings is rainbow by Claim 1. Due to the colors corresponding to $Y'_{1,1}$, one of these colorings is $W'_{1,1}$-canonical, so by Claim 2  none of the colorings is $X'_{1,1}$-canonical. Consequently, Claims 3-5 imply that
there is an integer $p_1$,  and a set $\{w_{1,1}(y)\ldots, w_{1, p_1}(y)\}$ for each $y \in Y'_{1,1}$, whose vertices correspond to the $W'_{1,1}$-canonical colors of $w_1(y)$.  Moreover, $w_{1,1}(y)=y$,
\begin{equation}\label{01192}
\mbox{$w_{1,j}(y) \ne w_{1,j}(y')$ for $y \ne y'\quad$ and
 $\quad w_{1,j}(y) \ne w_{1,j'}(y)$ for $j \ne j'$.}
\end{equation}

  In other words the $(j+1)$st canonical color class contains all edges of the form $xw_1(y)w_{1,j}(y)$ for $x \in X'_{1,1}$ and $y \in Y'_{1,1}$.
 Let $|X'_{1,1}|= |Y'_{1,1}|=s'_{1,1}$ and $Y'_{1,1}=\{y_1,\ldots,y_{s'_{1,1}}\}$. Add the colors of monochromatic colorings to $Q$.

  {\bf Type-1 cleaning:} Recall that $w_{1,1}(y_h)=y_h$ for all $1\leq h\leq s'_{1,1}$.
  By Lemma~\ref{trimlemma} with $A_h=\{y_h,w_1(y_h)\}$ and $a_h=w_{1,2}(y_h)$,
we can renumber $y_h$ so that the sets $$\{y_1, w_{1}(y_1),w_{1,2}(y_1)\},\ldots,\{y_{s'_{1,1}/3}, w_{1}(y_{s'_{1,1}/3}),w_{1,2}(y_{s'_{1,1}/3})\}$$
are pairwise disjoint. Applying  Lemma~\ref{trimlemma}  with $A_h=\{y_h,w_1(y_h),w_{1,2}(y_h)\}$ and $a_h=w_{1,3}(y_h)$,
then with $A_h=\{y_h,w_1(y_h),w_{1,2}(y_h),w_{1,3}(y_h)\}$ and $a_h=w_{1,4}(y_h)$, and so on, we obtain that
for $s''_{1,1}=\left\lceil\frac{s'_{1,1}}{3^{p_1}}\right\rceil$
we can renumber $y_h$ so that the sets $$\{y_1, w_{1}(y_1),w_{1,2}(y_1),\ldots,w_{1,p_1}(y_1)\},\ldots,
\{y_{s''_{1,1}}, w_{1}(y_{s''_{1,1}}),w_{1,2}(y_{s''_{1,1}}),\ldots, w_{1,p_1}(y_{s''_{1,1}})  \}$$
are pairwise disjoint. Let
    $Y''_{1,1}=\{y_1,\ldots,y_{s''_{1,1}}\}    $ and $X''_{1,1}$ be any subset
   of $X'_{1,1}$ of size $s''_{1,1}$.
   This is the {\em type-1 cleaning}.

{\bf  Type-2 cleaning:} Note that we automatically have $w_{1,j}(y) \cap X''_{1,1}=\emptyset$,
 since $xw_1(y)w_{1,j}(y) \in H_1$ for all $x \in X''_{1,1}$ so in particular, these three vertices are distinct.
Since for every $1\leq j\leq p_1$ all vertices $w_{1,j}(y_h)$ are distinct,
at most $|Q|\leq \ell$ of them are in $Q$. Deleting from $Y''_{1,1}$ the at most $p_1|Q|$ vertices $y_h$ such that
$$\{y_1, w_{1}(y_1),w_{1,2}(y_1),\ldots,w_{1,p_1}(y_1)\}\cap Q\neq \emptyset,$$
we  obtain a $Y_{1,1}\subset Y''_{1,1}$ such that for distinct $y\in Y_{1,1}$ the sets
$\{y, w_{1}(y),w_{1,2}(y),\ldots,w_{1,p_1}(y)\}$ are disjoint from each other and  from $Q$ and $X''_{1,1}$.
 Then we  choose any  $X_{1,1}\subset X''_{1,1}$ with $|X_{1,1}|=|Y_{1,1}|$.
This is the {\em type-2 cleaning}.

Now define $G_{1,2}$ to be the complete bipartite graph with parts $X_{1,1}$ and $W_2$ and repeat the
cleaning procedures above to obtain the integer $p_2$,  subsets
$X_{1,2}  \subset X_{1,1}$ and $Y_{1,2} \subset Y_{1,1}$ and vertices $w_{2,j}(y)$ that are distinct for distinct $y$ and also distinct from $w_{1,j'}(y')$ if $y \ne y'$.
Continuing in this way we obtain sets $X_{1,1} \supset X_{1,2} \supset \cdots \supset X_{1,p}:=X_2$ and
$Y_{1,1} \supset Y_{1,2} \supset \cdots \supset Y_{1,p}:=Y_2$,
 $\bV_2=\{(i,j_i): i\in [p], j_i \in [p_i]\}\subset [\ell]^2$ and vertices $w_{\bv}(y)$ for $\bv \in \bV_2$ and $y \in Y_2$ with $w_\bv(y) \not\in \{w_\bv(y'), w_i(y')\}$ for $y \ne y'$.

Given a vector ${\bf x}$ let ${\bf x}*j$ be the vector obtained from ${\bf x}$ by adding a new last coordinate with entry $j$ (for example if ${\bf x}= (3,7)$ then ${\bf x}*4=(3,7,4)$).
For $\bv \in \bV_2$, set $W_\bv=\cup_{y \in Y_2} w_\bv(y)$.
 Let us
 also construct the auxiliary digraph $D_2$ with vertex
set
$Y_2\cup  \bigcup_{i=1}^p W_i \cup \bigcup_{\bv \in \bV_2}W_\bv$
 with edges of the form $y w_i(y)$ for all $y,i$ and  $w_i(y)w_{i,j}(y)$ for $i \in [p]$ and $j \in [p_i]$.
 Because of cleanings, $D_2$ is  the vertex disjoint union of  homomorphic images of trees of height at most two, one for each $y\in Y_2$.


 \begin{claim}\label{cla3'} $|Q|\leq k$.
\end{claim}
{\em Proof:} By the definition of monochromatic colorings and by construction, for each $x\in X_2$ and each $w\in Q$,
$xw\in \partial H_1$ and the codegree of $xw$ is larger than $3k$. So we simple embed $T$ into the complete bipartite graph with partite sets
$X_2$ and $Q$, and then expand it.\qed

 To summarize, we have a set of (one dimensional) vectors $\bV_1=\{(1), \ldots, (p)\}$, nonnegative integers $p_\bv\le \ell$ for each $\bv \in \bV_1$ and

\begin{center}
\begin{tabular}{lp{5in}}
$\bullet$ & $\bV_2 =\cup_{\bv \in \bV_1}\{\bv*i: i\in [p_\bv]\} \subset [\ell+1]^2$, \\
$\bullet$ & $X_2 \subset X$ and $Y_2 \subset Y$, \\
$\bullet$ & vertices $w_\bv(y)$ with $w_{\bv}(y) \ne w_{\bv'}(y')$ if $y \ne y'$ and $\bv, \bv' \in \bV_1 \cup \bV_2$, \\
$\bullet$ & edges $xw_\bv(y) w_{\bv*i}(y) \in H_1$ for all $x \in X_2, y \in Y_2, \bv \in \bV_1, i \in [p_\bv]$ (so $\bv*i \in \bV_2$), \\
$\bullet$ & a digraph $D_2$ with vertex set
$Y_2 \cup \bigcup_{y \in Y_2, \bv \in \bV_1 \cup \bV_2} w_\bv(y)$
 and edges $y w_{\bv}(y)$ for $y \in Y_2, \bv \in \bV_1$ and $w_\bv(y) w_{\bv'}(y)$ as long as $\bv'=\bv*j$ for some $j \in [p_{\bv}]$, \\
$\bullet$ & the set $Q$ of all ``central'' vertices in monochromatic colorings, and $|Q|\leq k$.
\end{tabular}
\end{center}

\medskip

{\bf General Setup: }  Let $t \le \ell_2+1$ and suppose we have the following:

\begin{center}
\begin{tabular}{lp{5in}}
$\bullet$ & $\bV_t \subset [\ell+1]^t$, \\
$\bullet$ & $X_t \subset X$ and $Y_t \subset Y$, \\
$\bullet$ & for all $\bv \in \cup_{i=1}^t\bV_t$ and $y \in Y_t$ a vertex $w_{\bv}(y)$ such that for $y\ne y'$, $w_\bv(y)\ne w_{\bv'}(y')$, \\
$\bullet$ & edges $xw_{\bv}(y)w_{\bv*i}(y) \in H_1$ for all $x \in X_t, y \in
Y_t, \bv \in \cup_{j=1}^{t-1}\bV_j, i \in [p_\bv]$, \\
$\bullet$ & a digraph $D_t$ with vertex set
$$Y_t \cup \bigcup_{y \in Y_t, \bv \in \cup_{i=1}^t \bV_t} w_\bv(y)$$
 and edges $w_\bv(y) w_{\bv'}(y)$ as long as $\bv'=\bv*j$ for some $j \in [p_{\bv}]$ (define $y:=w_{{\bf \emptyset}}(y)$), \\
$\bullet$ & the set $Q$ of all "central" vertices in monochromatic colorings, and $|Q|\leq k$.
\end{tabular}
\end{center}

We will now show how to construct the same setup with $t+1$.
\medskip

Let  $\bV_t=\{\bv(1), \ldots, \bv(m(t))\}$.
Consider the complete bipartite subgraph $G_{t,1}$ of $\partial H_1$ with parts $X_t$ and $W_{\bv(1)}=\{w_{\bv(1)}(y): y \in Y_t\}$.
We apply the basic cleaning procedure to $G_{t,1}$ and obtain subsets $X'_{t,1} \subset X_t$ and $Y'_{t,1} \subset  Y_t$ and
 colorings $M_1, \ldots, M_{\ell+1}$ of the edges of $G_{t,1}$
  that are rainbow, canonical, or monochromatic. By Claim 1, no coloring is rainbow. By construction, we already have one
$W_{\bv(1)}$-canonical coloring obtained by considering the in-neighbors of $w_{\bv(1)}(y)$ in $D_t$.
By Claim 4, $R\neq\emptyset$ and thus $\ell_1\leq \ell$. We may assume that $M_1$ is $W_{\bv(1)}$-canonical.
Hence by Claim 2 no $M_i$ is $X'_{t,1}$-canonical. By Claim 3, the number of monochromatic colorings is at most $\ell_1-1\leq \ell-1$,
which means that the number of  $W_{\bv(1)}$-canonical colorings is at least $(\ell+1)-(\ell-1)=2$.
Consequently, there is a positive integer $p_{\bv(1)}$
 and $p_{\bv(1)}$ colorings (excluding the $W_{\bv(1)}$-canonical coloring given by the in-neighbors of $W_{\bv(1)}$) that are
 $W_{\bv(1)}$-canonical and the remaining $\ell+1-(p_{\bv(1)}+1)$ colorings are monochromatic.
We also have vertices $w_{\bv(1) * i}(y)$ for all $i=1, \ldots, p_{\bv(1)}$ which are distinct for distinct $y$ and distinct $i$.
As  before, for each $j \in [p_{\bv(1)}]$,  the $j$th canonical color class consists of all (hyper)edges of the form
$x w_{\bv(1)}(y) w_{\bv(1)*j}(y)$ for all $ x \in X'_{t,1}, y \in Y'_{t,1}$.

Next we perform the type-1 cleaning procedure (using Lemma~\ref{trimlemma}) to make sure that if $y \ne y'$ then $w_{\bv(1)*i}(y) \not\in \{w_\bv(y'), y'\}$
for any $\bv \in \bV_1 \cup \ldots \cup \bV_{t}$.
 This results in subsets $X''_{t,1} \subset X'_{t,1}$ and $Y''_{t,1} \subset Y'_{t,1}$. Finally, we perform the type-2 cleaning procedure to obtain $X_{t,1} \subset X''_{t,1}$ and $Y_{t,1} \subset Y''_{t,1}$ so that these sets do not contain any vertices that correspond to monochromatic colorings in any previous round.
 Add the central vertices of the monochromatic colorings to $Q$. Repeating the proof of Claim~6, we still have $|Q|\leq k$.

Now we repeat these procedure with $\bv(2)$ to obtain $X_{t,2} \subset X_{t,1}$ and $Y_{t,2}  \subset Y_{t,1}$.  Finally we perform this procedure with $\bv(m(t))$ to obtain
$X_{t+1}=X_{t, m(t)}$ and $Y_{t+1}=Y_{t, m(t)}$ and
$\bV_{t+1}=\cup_{\bv \in \bV_t}\{\bv*i:  i\in [p_\bv]\}$. We also have vertices $w_\bv(y)$ for every $y \in Y_{t+1}$ and $\bv \in \bV_{t+1}$ that are distinct for distinct $y$ and a digraph $D_{t+1}$ defined in the obvious way which consists of the vertex disjoint union of   homomorphic image of trees of height $t+1$,
one for each $y \in Y_{t+1}$.  Edges of the digraph encode the canonical colorings, as in the case $t=1,2$.

We repeat this procedure till we obtain sets $X_{g},Y_{g}, D_g$, for $g:=\ell_2+1$. By Claim 5, the outdegree of vertex $w_{\bv}(y) \in V(D_g)$ is
$$p_\bv \ge (\ell+2-\lambda)-1=\ell+1-\lambda.$$
Note that this is one less than the bound for $p$ because we have one in-neighbor that accounts for one canonical coloring.

\subsection{Embedding $T^+$ using $D_g$}
In this section we use the properties of $D_g$ to embed $T^+$ in $H_1$.
Our plan is to place the edges of $R$
on the edges of $D_g$ and to place the vertices of $I$ onto some vertices in $X_{g}$.
Let $T_1=T-L$.  Consider every tree  in the forest $R$ as a (directed) rooted tree $R_i$ with root $r_i$ which is a vertex in $V(R_i)$
 of the largest degree in $T_1$. Suppose we have $h$ such trees. By Lemma~\ref{treelem},
$$ d_{T_1}(r_i)\leq  \ell-\lambda 
  \quad\mbox{for all}\quad 1\leq i\leq h.$$


For each $ y \in Y_g$, let $D_g(y)$ be the component of $D_g$ containing $y$.
Choose $h$ vertices $y_1,\ldots,y_h\in Y_g$ arbitrarily,  and for  $1\leq i\leq h$ we will embed $R_i$ into ${D}_g(y_i)-w_1(y_i)$
 as follows (we exclude $w_1(y_i)$ because we will use $w_1(y_i)$ later in the embedding of $T^+$).
 Place $r_i$ on $y_{i}$. Suppose $r_i$ has $u$ out-neighbors in $R_i$. By construction,
 $y_{i}$ has $p\geq\ell+2-\lambda$ outneighbors in $D_{g}(y_i)$. So by~\eqref{root}, we can place
 the outneighbors in $R_i$ of $r_i$ on  outneighbors of $y_{i}$ in $D_{g}(y_i)$. 
 Then we place the outneighbors of placed vertices and so on. The general situation is that
some $v\in V(R_i)$ is placed on some $w_{\bv}(y)$ and has $u$ outneighbors in $R_i$.
By Lemma~\ref{jl1}, $\ell\geq 2\ell_2$. By Lemma \ref{lambda-lemma}, $\lambda = \lambda(R) \leq |R|/2 = \ell_2/2$.
So $w_\bv(y)$ has  $p_\bv \ge\ell+1-\lambda\geq \frac{3\ell_2}{2}+1$ outneighbors in $D_{g}(y)$.
At most $\ell_2-u$ of them are already occupied by previously embedded vertices. This leaves more than
$u$ available outneighbors of $w_\bv(y)$ to place the outneighbors in $R_i$ of $v$ on them.

\medskip
After placing all vertices in $V(R)$, we call a vertex of $H_1$ {\em free}, if it is not occupied by vertices
in $V(R)$ and is not the outneighbor of any occupied vertex in $D_g$. By construction,
there are at most $|V(R)|\ell_2\ell\leq \ell^3$ non-free vertices.
We now
place the vertices of $I$ on arbitrary distinct vertices in $X_g$
(they are all free at this moment by construction). Then we place the vertices of $D$ on distinct free
vertices in $Y_g$. Let $\phi$ be the embedding we are producing. We will assume below that each $a\in  I$ was placed on $\phi(a)\in Y_g$.
This yields an embedding of $T_1$ into $\partial H_1$.
In what follows, say that a pair $xy$ is {\em expanded} to a triple $xyz$.
 Our next goal will be to expand the
edges of $T_1$. After that we will embed the edges of $T-T_1$ and  expand them (these are the edges incident to $L$).

\medskip
Since the codegree of every edge in $D_g$ is at least $|X_g|$, we do not worry about
expanding the edges in $R$: we can do it greedily at the end. Recall that vertices in $D$ are adjacent only to $I$.
We need to expand the $|I|+|D|+h-1$ edges connecting $I$ with $D \cup V(R)$. For every
host $y$ of a vertex  $a\in D$ and the host $x$ of one neighbor $a'$ of $a$ in $I$,
we expand the edge $yx$ to $\{x,y,w_1(y)\}$. So the number of edges of $T_1-R$ not yet expanded is
$|I|+h-1=\ell_1+h-1$. Since the sets $V({D}_g(y))$ are disjoint for distinct $y\in Y_g$, expanding the edges
incident with $\phi(a)$ for $a\in D$ is easy: we simply use the vertices $w_2(\phi(a)),w_3(\phi(a))$ and so on.
Since the number of such edges is at most $\ell_1-1\leq p-2$, no problem arises.

When we expand an edge $yx$ where $x$ is the host of some $a\in I$ and $y$ is the host of some
$b\in V(R_i)\subset V(R)$, we need some more care, since some outneighbors of $y$ in $D_g(y)$
can be occupied. For $i=1,\ldots,h$ let $U(i)=|R_i|+|E_{T_1}(I,V(R_i))|$. Then
\begin{equation}\label{ui}
 \sum_{i=1}^hU(i)=|T_1-A_0|\leq |I|+|R|+h-1=\ell+h.
 \end{equation}
Order the $R_i$s so that $U(1)\geq U(2)\geq \ldots\geq U(h)$ and expand the edges incident to $R_i$s in the
reverse order.
Since each $b\in V(R)$ is adjacent to some $a\in I$, $U(i)\geq 3$ for every $i$.
Suppose that it is now the turn to expand the edges incident to $R_i$ and $i\geq 2$.
Then
$U(i)\leq U(2)\leq \frac{\ell+h-3(h-2)}{2}\leq \frac{\ell+2}{2}$.
We expand the edges one by one. Suppose we need now to expand $w_\bv(y)x$, where
$w_\bv(y)$ is the host  of a vertex $b\in V(R_i)$ (possibly $\bv=\emptyset$ in which case by convention $w_\bv(y)=y$ and $p_\bv=p$).
The outdegree in $D_{g}$ of $w_\bv(y)$ is
 $p_\bv\ge \ell+1-\lambda$. At most $|R_i|$ of the outneighbors of $w_\bv(y)$ are occupied.
If we already expanded some edges incident with $R_i$, they block at most $|E_{T_1}(I,V(R_i))|-1$
outneighbors of $y$. Consequently, we have at least
$$(\ell+1-\lambda)-(U(i)-1)\geq \frac{\ell}{2}-\lambda+1\ge
\frac{\ell}{2}-\frac{\ell_2}{2}+1>0$$ free outneighbors of $y$, and any free outneighbor may be used to expand   $w_\bv(y)x$.

Finally, we work with $R_1$. It is possible that $U(1)$ is as large as $\ell_2+\ell_1$.
On the other hand, we have never used the universal vertices for monochromatic multicolorings,
and this is the time to use them.
Now for each $a\in V(R_1)$ and $x\in X_g$, the pair $\phi(a)x$ has
  $1+\ell=\ell_1+\ell_2$ different colors in the canonical multicoloring (including any universal vertices), which means $1+\ell$ possibilities to expand
$\phi(a)x$.  Since the number of edges of $U(1)$ to be embedded when we embed $R_1$ is at most $\ell+1$, we can perform the embedding greedily.

Having embedded and expanded $T_1$, we work with $L$. Since $Y_g$ is large, one by one, take $c\in L$, place it on a free $y\in Y_g$ and expand
the obtained edge $yx$ via $w_1(y)$.  \qed

\section{Proof of Theorem \ref{main2}}\label{ktt}

Suppose $\sigma(G^+) = 2$ and $|V(G)|=k$. Since the  $n$-vertex triple
system of all edges containing a fixed vertex does not contain $G^+$ with $\sigma(G^+) = 2$ (by definition),
$\mbox{ex}(n,G^+) \geq {n-1 \choose 2}$. Also if $\sigma(G^+) = 2$,
then either some vertex of
$G$ covers all but  one edge in $G$ (and this edge connects two leaves) or two non-adjacent
vertices of $G$ cover all edges of $G$. In the former case, $G$ is contained in the star-plus-one-edge graph $S^*_{k-1}$  and in the latter, $G$ is
contained in $K_{2,k-2}$. Thus it is enough to consider the cases $G=K_{2,k-2}$ and $G=S^*_{k-1}$.

Suppose we have an $n$-vertex $3$-graph $H$ not containing $G^+$ for $G\in \{S^*_{k-1},K_{2,k-2}\}$
 with $|H| = (1+\varepsilon){n\choose 2}$
where $\varepsilon > 0$ and $n$ is sufficiently large. It is enough to assume $k\geq 5$.
 Let $H'$ be obtained from $H$ by consecutive deletion of edges having a pair
of codegree one, so that the minimum codegree of edges in $H'$ is at least two.
If we deleted $m$ edges, then $|\partial H'| \leq {n\choose 2} - m$. Let $E$ be the set of edges of $H'$ in
which the codegrees of all pairs (in $H'$) are at most $3$ or at least two  pairs have codegree (in $H'$) exactly two. We claim that
\begin{equation}\label{P1}
|E|\leq |\partial H'|.
\end{equation}
To see this, define $\omega = \sum_{e \in H'} \sum_{f \subset e} 1/d(f)$, where $d(f)$ is the codegree of $f$ in $H'$.
By definition of $E$, for every $e \in E$ we have $\sum_{f \subset e} \frac{1}{d(f)} \geq 1$. Since $E \subset H'$,
we get $\omega \geq |E|$. By interchanging the sums, we see $\omega = |\partial H'|$:
\[ \omega = \sum_{f \in \partial H'} \sum_{e \supset f} \frac{1}{d(f)} = \sum_{f \in \partial H'} 1 = |\partial H'|.\]
Therefore $|E| \leq |\partial H'|$ as claimed.


Let $H'' = H' \backslash E$.
By (\ref{P1}),  $|H''|\geq \varepsilon {n \choose 2}$. By the definition of $E$, if $e\in H''$ and the codegrees in $H'$ of the vertex pairs
in $e$ are $c_1\leq c_2\leq c_3$, then $c_1\geq 2$, $c_2\geq 3$ and $c_3\geq 4$. Then
\begin{equation}\label{0122}
\mbox{
for each $e\in H''$, there is a triangle $T_e$
in $H'$ whose every edge shares $2$ vertices with $e$.}
\end{equation}
We partition $H''$ into three triple systems. Let
$H_1$ be the set of  $e\in H''$ containing a pair $f=f_e \subset e$ with $3 \le d_{H'}(f) \le 3k$,
$H_2$ be the set of  $e\in H''$ with one pair $f_e \subset e$ having $d_{H'}(f_e)=2$  and two pairs of codegree (in $H'$) at least $3k+1$, and
$H_3= \{e \in H'' : \delta_{H'}(e) \geq 3k+1\}$. By the definition of $H''$ we have $H_1\cup H_2\cup H_3=H''$ and one
of the three cases below must hold.

{\bf Case 1:}  $|H_1|\geq \frac{\varepsilon n^2}{9} $. Let $F=\{f_e: e\in H_1\}$ ($f_e$ is defined above) so that $|F|\ge |H_1|/3k\ge \varepsilon n^2/27k$.
For every $f\in F$, choose  $S_f \subset N_{H'}(f)$ with $|S_f|= 3$ such that $S_f\cap N_{H_1}(f)\neq\emptyset$ (we can do it, since by definition,
each $f\in F$ is $f_e$ for some $e\in H_1$).
By Lemma~\ref{kttlemma} applied to $F$, for a large $t$ there exists $K\subseteq F$ such that $K\cong K_{t,t}$
and for every $f\in K$,
  $S_f \cap V(K)=\emptyset$.
  By Theorem \ref{canon}, if $t$ is large enough, there exists $K'\cong K_{2k,2k} \subset K$
  and three disjoint list-edge-colorings $\chi_i : K' \rightarrow L_{K'}$ such that each $\chi_i$ is monochromatic or canonical, or some $\chi_i$ is rainbow.
Let $X = \{x_1,x_2,\dots,x_{2k}\}$ and $Y = \{y_1,y_2,\dots,y_{2k}\}$ be the parts of $K'$. If say
 coloring $\chi_1$ is rainbow, then clearly $K_{2,k - 2}^{+} \subset K_{2k,2k}^{+} \subset H_1$ and we are done when $G=K_{2,k - 2}$.
 Suppose $G=S^*_{k-1}$. By the construction of $K$, there is an edge $zx_1y_1\in H_1$
  such that $z\notin  V(K)$. By~\eqref{0122}, $H'$ contains a triangle $\{x_1y_1u_1,x_1zu_2,y_1zu_3\}$.
   For at most four values of $2\leq i\leq 2k$,  $\{z,u_1,u_2,u_3\}\cap \{y_i,\chi_1(x_1y_i)\}\neq\emptyset$.
 So, $H_1$ contains $(S^*_{2k-5})^+$ with the center $x_1$. Since $k\geq 5$, we are done.

 Suppose now that no coloring is rainbow. We have three possibilities.

\medskip
{\bf Case 1.1.} {\em $G=S^*_{k-1}$.} If some coloring $\chi_i$ is monochromatic, say,
$\chi_1(e) = \alpha$  for all $e \in K_{2k,2k}$,
then the edges $x_iy_i\alpha$  for $1 \leq i \leq k-1$ and the edge $x_1y_2\chi_2(x_1y_2)$  form a $(S^*_{k-1})^{+} \subset H'$ with
the center $\alpha$.
Otherwise, we may assume that $\chi_1$ is $X$-canonical. Let $\alpha_i$ be the color in $\chi_1$ common to every edge containing $x_i$.
Since $d_{H'}(y_1\alpha_1)\geq 2$, there is a vertex $w\neq x_1$ such that $wy_1\alpha_1\in H'$. By symmetry, we may assume that
$w\notin \{x_2,y_2,\ldots,x_k,y_k\}$.
Then the edges $y_1\alpha_i x_i$  for $2 \leq i \leq k-2$, $wy_1\alpha_1$, $x_1\alpha_1y_2$ and
$y_1x_1\chi_2(x_1y_1)$ form a $(S^*_{k-1})^{+} \subset H'$ with
the center $y_1$.

{\bf Case 1.2.} {\em $G=K_{2,k - 2}$ and some coloring $\chi_i$ is monochromatic}. If two or more of the colorings are monochromatic, say,
 $\chi_1(e) = \alpha$ and $\chi_2(e) = \beta$ for all $e \in K_{2k,2k}$ with $\alpha \neq \beta$,
then the edges $x_iy_i\alpha$ and $x_iy_{i + k}\beta$ for $1 \leq i \leq k$ form a $K_{2,k}^{+} \subset H'$.
If only one coloring is monochromatic, then the other two are canonical. We may assume $\chi_1(e) = \alpha$
for $e \in K_{2k,2k}$ and $\chi_2$ is $X$-canonical.
Let $\alpha_i$ be the color common in $\chi_2$ to every edge containing $x_i$. Then the edges $\alpha x_iy_i$ and
$\alpha_ix_iy_{k + 1}$ for $1 \leq i \leq k$ form a $K_{2,k}^{+}\subset H'$.


{\bf Case 1.3:} {\em $G=K_{2,k - 2}$ and no coloring is monochromatic.} This means all of the $\chi_i$ are canonical.
In particular, by symmetry, we can assume $\chi_1$ and $\chi_2$ are both $X$-canonical. If $\alpha_i$
is the common color of every edge on $x_i$ under $\chi_1$, and $\beta_i$ is the common color of every edge on $x_i$ under $\chi_2$,
then the edges $y_1x_i\alpha_i$ and $y_2x_i\beta_i$ for $1 \leq i \leq k$ form a $K_{2,k}^{+}  \subset H'$.
This finishes Case 1.

{\bf Case 2:}  $|H_2|\geq \frac{\varepsilon n^2}{9} $. By  the K\"{o}vari-S\'{o}s-Tur\'{a}n Theorem, for every $k$ there is $s(k)$ such that
every subgraph $M$ of $K_{s(k),s(k)}$ with at least $s(k)^2/2$ edges contains a $K_{2k,2k}$. Similarly to Case 1,
let $F=\{f_e: e\in H_2\}$, where $d_{H'}(f_e)=2$. For every $f\in F$, let  $S_f = N_{H'}(f)$.
By definition, $|S_f|= 2$ and $S_f\cap N_{H_2}(f)\neq\emptyset$.
Then $|F|\ge |H_1|/2\ge \varepsilon n^2/18$.
By Lemma~\ref{kttlemma} applied to $F$, for a large $t$ there exists $K\subseteq F$ such that $K\cong K_{t,t}$
and for every $f\in K$,  $S_f \cap V(K)=\emptyset$.
  By Theorem \ref{canon}, if $t$ is large enough, there exists $K_0\cong K_{s(k),s(k)} \subset K$
  and  disjoint list-edge-colorings $\chi_1$ and $\chi_2$ of  $K_0$  such that each $\chi_i$ is monochromatic or canonical, or some $\chi_i$ is rainbow.
Since each of the lists contains a color corresponding to an edge in $H_2$, we may assume that for at least of half of the edges  $f\in K_0$,
$f\cup \{\chi_1(f)\}\in H_2$. Then by the definition of $s(k)$, there exists $K'\cong K_{2k,2k} \subset K_0$ such that for every $f\in K'$,
$f\cup \{\chi_1(f)\}\in H_2$.
Now we repeat the proof of Case 1 word by word till (and including) Case 1.2, since in these subcases we we have used only two colorings.
In Case 1.3, the problem arises only when $\chi_1$ is $X$-canonical and $\chi_2$ is $Y$-canonical (or vise versa).
Let  
$\chi_2(x_1y_i) = \alpha_i$ for $1 \leq i \leq k$. Since $\chi_1$ is $X$-canonical, we have edges $y_ix_1\gamma \in H'$ for $1 \leq i \leq k$, where $\gamma$ is the common color of all edges on $x_1$ in {\color{red}$\chi_1$}. By construction, for every $1\leq i\leq k$, edge $x_1y_i\gamma$ is in $H_2$ and hence
 $d_{H'}(y_i\gamma) \geq 3k + 1$. Therefore we may choose vertices $\beta_1,\beta_2,\dots,\beta_k \in L_K \backslash \{y_1,y_2,\dots,y_k,\alpha_1,\alpha_2,\dots,\alpha_k\}$ such that $\gamma\beta_iy_i$ are all edges of $H'$. These edges together with
the edges $x_1y_i\alpha_i$ form $K_{2,k}^{+} \subset H'$.

{\bf Case 3:}  $|H_3|\geq \frac{\varepsilon n^2}{9} $.
If $|\partial H_3| > \frac{\varepsilon }{200k} n^2$, then similarly to Case 1,
for every $f\in \partial H_3$, choose  $S_f \subset N_{H'}(f)$ with $|S_f|= 3$ such that $S_f\cap N_{H_3}(f)\neq\emptyset$.
By Lemma \ref{kttlemma} applied to $F=\partial H_3$, for a large $t$ there exists $K\subseteq F$ such that $K\cong K_{t,t}$
and for every $f\in K$,  $S_f \cap V(K)=\emptyset$. From this point, we just repeat the proof of Case 1.

So $|\partial H_3| \leq \frac{\varepsilon }{200k} n^2$.
Then by Lemma \ref{fullsub}, $H_3$ contains an $8k$-full subgraph $H^*$ with at
least $|H_3| - 8k|\partial H_3|\geq \frac{\epsilon}{20}n^2$ edges.
Since $|\partial H^*| = o(n^2)$, $d_{H^*}(xy)\geq 2k$ for some  $xy\in\partial H^*$.
This means that the edge $xy$ in the graph  $\partial H^*$ is in at least $2k$ triangles.
So, $\partial H^*$ contains $S^*_k$ with the center $x$ and $K_{2,k}$ with the small partite set
$\{x,y\}$. This means $\partial H^*$ contains a copy of $G$.
Since $H^*$ is $8k$-full,
our copy of $G$ greedily extends to $G^+\subset H^*$. This finishes the main proof.

The jump in the Tur\'an number follows immediately by observing that if $\sigma(G^+)\ge 3$, then we may apply (\ref{triviallower})
and obtain ex$_3(n, G^+) \ge (2-o(1)){n \choose 2}$. \qed




\section{Concluding Remarks} \label{concluding}
$\bullet$ Our methods can be used to determine the order of magnitude of  the Tur\'an number of expansions of other bipartite graphs like the 3-dimensional cube and complete bipartite graphs.  These will be presented in a forthcoming paper.

$\bullet$
Our approach may also be suitable for other extremal problems on trees and forests in hypergraphs including the following conjecture of Kalai (see Frankl
and F\"{u}redi~\cite{FF}), extending the Erd\H{o}s-S\'{o}s Conjecture to $r$-graphs.  An $r$-tree is an $r$-graph with edges $e_1, \ldots, e_q$ where for each $i$, $e_{i} \cap (\cup_{j< i} e_j) \subset e_k$ for some $k<i$.

\begin{conjecture} {\bf (Erd\H os-S\'os for graphs and Kalai 1984 for $r \ge 3$)}\label{kalai}
Let $r \ge 2$ and $T$ be an $r$-tree on $v$ vertices.  Then
$$ex_r(n, T)\le \frac{v-r}{r} {n \choose r-1}.$$
\end{conjecture}

This conjecture has been solved for certain classes of trees (see~\cite{FF}).





\section{Acknowledgments}

We are grateful to Z. F\"{u}redi for informative discussions and helpful suggestions which improved this paper.

\end{document}